 \documentclass[11pt]{article}
 \usepackage[top=1in,bottom=1in,left=1.5in,right=1.5in]{geometry}
  \usepackage{amsmath,amssymb}
  \usepackage{latexsym}
\usepackage{epsfig}
\usepackage{epsf}
\usepackage{amsfonts}
\usepackage{graphicx}

  \newcommand{\F}{\mathbf{F}}

  \newcommand{\R}{\mathbb{R}}

  \renewcommand{\a}{\mathbf{a}}
  \renewcommand{\b}{\mathbf{b}}
  
  \newcommand{\e}{\mathbf{e}}
  \newcommand{\f}{\mathbf{f}}
  \newcommand{\g}{\mathbf{g}}

  \newcommand{\n}{\mathbf{n}}

  \newcommand{\p}{\mathbf{p}}

  \renewcommand{\u}{\mathbf{u}}
  \renewcommand{\v}{\mathbf{v}}
  
  \newcommand{\w}{\mathbf{w}}
  
  \newcommand{\x}{\mathbf{x}}
  
  \newcommand{\y}{\mathbf{y}}
  \newcommand{\z}{\mathbf{z}}
  \newcommand{\0}{\mathbf{0}}
  \newcommand{\1}{\mathbf{1}}
  
  \newcommand{\bet}{\boldsymbol{\beta}}
  \newcommand{\gam}{\boldsymbol{\gamma}}
  \newcommand{\zet}{\boldsymbol{\zeta}}
  \newcommand{\et}{\boldsymbol{\eta}}
  \newcommand{\xit}{\boldsymbol{\xi}}
  \newcommand{\Gam}{\mathbf{\Gamma}}

  \newcommand{\cA}{\mathcal{A}}
  \newcommand{\cB}{\mathcal{B}}

  \newcommand{\cF}{\mathcal{F}}
  \newcommand{\cG}{\mathcal{G}}

  \newcommand{\cP}{\mathcal{P}}

  \newcommand{\lan}{\langle}
  \newcommand{\ran}{\rangle}
  \newcommand{\an}[1]{\lan#1\ran}
  \def\diag{\mathop{{\rm diag}}\nolimits}
  \newcommand{\hs}{\hspace*{\parindent}}
  \newcommand{\proof}{\hs \textbf{Proof.\ }}

  \newcommand{\trans}{^\top}
  \newcommand{\qed}{\hspace*{\fill} $\Box$\\}

  \newcommand{\inn}{\mathrm{in}}

  \newcommand{\rH}{\mathrm{H}}

  \newcommand{\rS}{\mathrm{S}}

  \newtheorem{theo}{\bfseries \hs Theorem}[section]

  \newtheorem{prob}[theo]{\bfseries \hs Problem}
  \newtheorem{lemma}[theo]{\bfseries \hs Lemma}
  \newtheorem{algo}[theo]{\bfseries \hs Algorithm}
  \newtheorem{claim}[theo]{\bfseries \hs Claim}
  \newtheorem{corol}[theo]{\bfseries \hs Corollary}

  \numberwithin{equation}{section} 

\begin{document}

 \title{Maximizing Sum Rates in\\
 Gaussian Interference-limited Channels\thanks{Dedicated to the memory of Samuel Karlin}}

 \author{
 Shmuel Friedland\thanks{Department of Mathematics, Statistics and Computer Science,
   University of Illinois at Chicago, Chicago, Illinois 60607-7045,
   USA, and Berlin Mathematical School, Berlin, Germany,
   e-mail:friedlan@uic.edu} \, and \,
   Chee Wei Tan\thanks{Electrical Engineering Department, Princeton University, NJ
   08544, USA, e-mail:cheetan@princeton.edu}}

 \date{June 27, 2008}
 \maketitle

 \begin{abstract}
 We study the problem of maximizing sum rates in a Gaussian interference-limited channel that models
 multiuser communication in a CDMA wireless network or DSL cable binder.
 Using tools from nonnegative irreducible matrix theory, in particular the Perron-Frobenius Theorem
 and the Friedland-Karlin inequalities, we provide insights into the structural property of optimal
 power allocation strategies that maximize sum rates. Our approach is similar to the treatment of
 linear models in mathematical economies, where interference is viewed in the context of competition.
 We show that this maximum problem can be
 restated as a maximization problem of a convex function on a
 closed convex set.
 We suggest three algorithms to find the exact
 and approximate values of the optimal sum rates. In particular, our algorithms exploit the eigenspace
 of specially crafted nonnegative {\it interference matrices}, which, with the use of standard
 optimization tools, can provide useful upper bounds and feasible solutions to the nonconvex problem.



 \end{abstract}

 \section{Introduction}

Interference control is a major and fundamental issue in the design of a multiuser communication system, where interference is a central phenomenon. The goal of an interference-limited system is to develop a resource allocation scheme that is able to jointly optimize the performance of all users in the presence of mutual interference. To meet the growing demand of untethered applications, transmission rate adaptation by power control has emerged as active research issues in multiuser systems, which often impact the design of higher networking layers. This is due to innovations at the link layer and the need for utility maximization over networks \cite{Srikant03,ChiangLow07}, where a system-wide performance objective has to be maximized or the quality of service for each individual user has to be met.

Power control is often used in multiuser communication system such as cellular, ad-hoc networks and digital subscriber lines (DSL) to provide a high signal-to-noise ($\mathsf{SNR}$) ratio to each user for a reliable connection. A higher $\mathsf{SNR}$ allows a system with link adaptation to transmit at a higher data rate, thus leading to a greater spectral efficiency. However, in a multiuser system, the performance of each user depends not only on its power allocation, but also on the power allocation of other users. There is thus a tradeoff among the data rates of all users, when a single system-wide performance metric or network-wide utility has to be optimized. Intelligent power control is necessary to jointly maximize the overall performance and mitigate the interference of all users.

The motivation of this work comes from maximizing sum rate or data throughput in wireless communication channels or DSL. Due to the broadcast nature of radio transmission, data rates in a wireless network are affected by multiuser interference. This is particularly true in Code Division Multiple
 Access (CDMA) systems, where users transmit at the same time over the same frequency
 bands and their spreading codes are not perfectly orthogonal. Transmit power control
 is often used to control signal interference to maximize the total transmission rates
 of all users. In a DSL broadband access system, multiple users share a cable binder and communicate with a central office through a common spectrum \cite{Luo08,Yu02}. Signals bundled in the close proximity of a cable binder interfere with one another due to electromagnetic coupling between neighboring twisted-pairs. In this paper, we focus primarily on sum rate maximization in a wireless communication channel over a single frequency, but we will also show how our results can be easily extended to find application to optimal spectrum management over multiple frequencies as in a DSL channel.

From an information theoretic viewpoint, finding the capacity region and optimal code design of the interference channel is still an open problem even for the simple case of two users. A practical approach that leads to an inner bound of the unknown capacity region or achievable rate region is to adopt a noncooperative approach that treats interference as additive Gaussian noise, i.e., no multiuser detection or interference cancellation is used. The Gaussian assumption on the interference distribution is justified, especially when the number of interfering users is large, e.g., in existing DSL cable binder, there are around fifty to hundred twisted pairs in a cable binder.
Such an approach is also relevant in practical systems, because the signaling overhead involved in coordination is minimal (leading to simple channel access protocols without a centralized coordinator and the use of single user codes) and transceiver techniques with low complexity can be employed at each user.

The rate region characterizes all possible data rate combinations among all users subject to power constraints while making the decoding error probability arbitrary small. Obtaining the optimal transmit power control scheme that maximizes sum rates requires solving a nonconvex problem \cite{Tan05,Chiang07,Ebrahimi06,Luo06,Yu02,Charafeddine07}. The computational aspect of an exhaustive search over all possible power allocations is very high. As the number of users increases, it becomes hard to obtain optimal solutions in polynomial time. While global optimality in power control allocation is highly desired, it is often of necessity to find local optimal solution with performance guarantees in the presence of nonconvexities. However, weighted sum rate is not only difficult to optimize, but the problem may even be hard to approximate \cite{Luo08}.

 In this paper, using tools from nonnegative irreducible matrix theory, in particular the Perron-Frobenius Theorem and the Friedland-Karlin inequalities, we provide insights into the structural property of optimal power allocation strategies that maximize sum rates. Our approach is similar to the treatment of linear models in mathematical economies, where interference is viewed in the context of competition. In particular, we show that we can convert the problem of maximizing sum rate
 to a problem of maximizing a convex function on a closed
 unbounded convex set. This approach is based on construction of {\it interference matrices} that characterize the problem, which is similar to the role of input-output matrices in mathematical economies. More importantly, it leads to the design of new
 algorithms for finding the exact and approximate solutions to
 the maximum problem. By focusing on fast and simple algorithms, interference among multiple users in a common spectrum can be mitigated so that efficient data multiplexing can be realized.

 We now state briefly the maximal problem.  (See \S2 for all
 definitions, notations and motivations.)
 Let $F=[f_{ij}]_{i=j=1}^L,\v=(v_1,\ldots,v_L)\trans$ be an
 $L\times L$ matrix, with zero diagonal and positive off
 diagonal elements, and a positive vector, respectively.  Define
 the following transformation $\gam(\p)=(\gamma_1(\p),\ldots,\gamma_L(\p))\trans$
 from the set of nonnegative vectors $\p\in\R_+^L$ to itself:
 \[\gamma_l(\p)=\frac{p_l}{\sum_{m=1}^L f_{lm}p_m +v_l}, \textrm{ for
 } l=1,\ldots,L, \textrm{ where } \p=(p_1,\ldots,p_l)\ge \0.\]
 Let $\bar\p=(\bar p_1,\ldots,\bar p_L)\trans$ be a given
 positive vector.  For a given probability vector
 $\w=(w_1,\ldots,w_L)\trans$, let
 \begin{equation}\label{defPhiw}
 \Phi_{\w}(\gam):=\sum_{i=1}^L w_i\log(1+\gamma_i), \textrm{ where }
 \w,\gam\in\R_+^L.
 \end{equation}
 Then the sum rate maximization problem in Gaussian interference-limited channels is
 \begin{equation}\label{defopt}
 \max_{\0\le p_l \le\bar{p}_l \; \forall \, l} \Phi(\gam(\p)).
 \end{equation}
The exact solution to this problem is known to be (strongly) NP-hard \cite{Luo08}. An often used classical technique to tackle nonconvexity is the Lagrange duality. The shortcoming of using this approach is however that there can exist a duality gap between the global optimal primal and dual solution of (\ref{defopt}) \cite{Luo08}. Also, finding an optimal primal solution given an optimal dual solution, or vice versa, is in general very difficult. In contrast to simply resorting to classical global optimization techniques, we will highlight in this paper the insightful role of using irreducible nonnegative matrix theory to tackle nonconvexity of the form given by (\ref{defopt}).


 For any vector $\tilde\gam=(\tilde\gamma_1,\ldots,\tilde\gamma_L)\trans\in\R^L$ let
 $e^{\tilde\gam}=(e^{\tilde\gamma_1},\ldots,e^{\tilde\gamma_L})\trans$.
 The purpose of this paper to show that the sum rate maximization problem
 equivalent to the maximization problem of the convex
 function $\Phi_{\w}(e^{\tilde\gam})$, where $\tilde\gam$
 varies on a closed unbounded convex domain $D(\{F\}\subset\R^L)$.
 (See (\ref{defconvset}).)  Note that this formulation shows
 that we can not expect, in general, that a maximal solution $\p^{\star}$
 is unique.
 We show how to approximate this domain by a convex polytope
 $D(\zet_1,\ldots,\zet_M,K)$, formed mainly by the supporting hyperplanes
 of $D(\{F\})$.
 This approximation is obtained by using the Friedland-Karlin inequality
 stated in \cite{Friedland75}.
 Hence, it is easier to find the
 maximum $ \Phi(e^{\tilde\gam})$ on $D(\zet_1,\ldots,\zet_M,K)$,
 than on $D(\{F\})$.  Finally we approximate the function
 $ \Phi(e^{\tilde\gam})$ by the linear function $\w\trans
 \tilde\gam=\sum_{l=1}^L w_l\tilde\gam_l$.  Then the maximum
 of $\w\trans\tilde\gam$ on $D(\zet_1,\ldots,\zet_M,K)$ is a
 linear programming problem, which can be solved in polynomial
 time using the ellipsoid algorithm \cite{Kha}.

 We survey briefly the contents of the paper.
 In \S2 we state definitions, notations and a short motivation.
 We give a characterization of
 the image of the multidimensional box $[\0,\bar\p]\subset
 \R^L_+$ by the map $\gam$, in terms of certain inequalities of
 the spectral radii of corresponding nonnegative matrices.
 In \S3 we study the sum rate maximization problem in wireless channel.
 We give necessary and sufficient conditions for an extremal
 point $\p\in [\0,\bar\p]$ of $\Phi_{\w}(\gam(\p))$ to be a
 local maximum.  In \S4 we study several relaxation versions of
 the sum rate maximization problem.  These relaxation version
 are simpler than the original maximal problem.
 When replacing the maximal function $\Phi_{\w}$ by the linear
 function $\w\trans\tilde\gam$, we have a closed-form solution under certain
 conditions on $\w$. This solution is obtained by using the Friedland-Karlin inequality
 stated in \cite{Friedland75}.
 In \S5 we give three algorithms
 to solve the sum rate maximization problem and its two other
 approximation versions described above. In \S6 we apply our results to sum rate maximization in DSL channel. We conclude the paper in \S7.
 In \S A, viewed as an appendix, we restate some useful results for \cite{Friedland75}
 and give several applications and extensions, which are needed
 in this paper.

 \section{Notations and preliminary results}
\label{modelsection}
 As usual, let $\R^{m\times n} \supset\R_+^{m\times n}$ denote the set of $m\times n$
 matrices and its subset of nonnegative matrices.  For
 $A,B\in\R^{m\times n}$, we denote $A\le B$ if $B-A\in
 \R_+^{m\times n}$. We denote $A \lneq B,A<B$ if $B-A$ is a nonzero nonnegative and positive
 matrix, respectively.  We denote the entries of a matrix $A\in \R^{m\times n}$ by the small letters, i.e $A=[a_{ij}]_{i,j=1}^{m,n}$.
 Identify $\R^{m}=\R^{m\times 1},
 \R_+^{m}=\R_+^{m\times 1}$.

 A column vector is denoted by the
 \textbf{bold} letter $\x=(x_1,\ldots,x_m)\trans \in \R^m$.
 Then $e^{\x}:=(e^{x_1},\ldots,e^{x_m})\trans$.
 For $\x>\0$, we let
 $\x^{-1}:=(\frac{1}{x_1},\ldots,\frac{1}{x_m})\trans$ and
 $\log\x=(\log x_1,\ldots,\log x_L)\trans$.
 For any $\x=(x_1,\ldots,x_L)\trans\in\R^L$, we let
 $e^{\x}=(e^{x_1},\ldots,e^{x_L})\trans$.
 Let
 $\mathbf{x} \circ \mathbf{y}$ denote the Schur product of the vectors
 $\mathbf{x}$ and $\mathbf{y}$, i.e., $\mathbf{x} \circ \mathbf{y}=[x_1 y_1, \dots, x_L y_L]^T$.
 Let $\1=(1,\ldots,1)\trans\in \R^{L}$.
 For $\underline{\p}\le\bar\p\in \R^L$, denote by
 $[\underline{\p},\bar\p]$ the set of all $\x\in\R^L$
 satisfying $\underline{\p}\le\x\le\bar\p$.
 For a vector $\y=(y_1,\ldots,y_L)\trans$, denote by
 $\diag (\y)$
 the diagonal matrix $\diag(y_1,\ldots,y_L)$.
 We also let $(B\mathbf{y})_l$ denote the $l$th element of $B\mathbf{y}$.
 The Perron-Frobenius eigenvalue of a nonnegative matrix $F$ is denoted as $\rho(F)$,
 and the Perron (right) and left eigenvector of $F$ associated with
 $\rho(F)$ are denoted by $\mathbf{x}(F)$ and $\mathbf{y}(F)$
 (or simply $\mathbf{x}$ and $\mathbf{y}$ when the context is clear) respectively.
 Assume that $F$ is a nonnegative irreducible matrix.  Then $\rho(F)$ is simple and positive, and
 $\x( F),\y(F)>\0$.  We will assume the normalization:
 $\x(F)\circ\y(F)$ is a probability vector.
 The super-script $(\cdot)^{\trans}$ denotes transpose.
 For a positive integer $n$, denote by $\an{n}$ the set
 $\{1,\ldots,n\}$.  Let $P: X\to Y$ be a mapping from the space $X$ to the
 space $Y$.  For a subset $Z\subset X$, we denote by $P(Z)$
 the image of the set $Z$.

 Consider a wireless interference channel, e.g., cellular network, with $L$ logical
 transmitter/receiver pairs. In the cellular uplink case, all logical receivers may reside in the same
 physical receiver, i.e., the base station. The interference-limited channel with $L$ users can be modeled as a Gaussian interference channel given by the baseband signal model:
 \begin{equation}
 \label{channelmodel}
 \boldsymbol{y}_l = \boldsymbol{h}_{ll} \boldsymbol{x}_l + \sum_{j \ne l} \boldsymbol{h}_{lj} \boldsymbol{x}_j + \boldsymbol{z}_{l},
 \end{equation}
where $\boldsymbol{y}_l \in \mathbb{C}^{1 \times 1}$ is the received signal of the $l$th user, $\boldsymbol{h}_{lj} \in \mathbb{C}^{1 \times 1}$ is the channel coefficient between the transmitter of the $j$th user and the receiver of the $l$th user, $\boldsymbol{x} \in \mathbb{C}^{N \times 1}$ is the transmitted (information carrying) signal vector, and $\boldsymbol{z}_l$'s are the i.i.d. additive complex Gaussian noise coefficient with variance $n_l/2$ on each of its real and imaginary components. The first term on the right-hand side of (\ref{channelmodel}) represents the desired signal, whereas the second term represents the interfering signals from other users. At each transmitter, the signal is constrained by an average power constraint, i.e., $\mathbb{E}[|\boldsymbol{x}_l|^2]=p_l$, which we assume to be upper bounded by $\bar{p}_l$ for all $l$.

 The vector $(p_1, \dots, p_L)^{\trans}$ is the transmit power vector and is the parameter of interest in this paper. Let $G=[g_{lj}]_{l,j=1}^L>0_{L\times L}$
 represents the channel gain, where $g_{lj}=|\boldsymbol{h}_{lj}|^2$ is the channel gain from the $j$th
 transmitter to the $l$th receiver, and $\n=(n_1,\ldots,n_L)\trans>\0$, where $n_l$
 is the noise power at the $l$th receiver. Assuming a linear matched-filter receiver, the Signal-to-Interference
 Ratio ($\mathsf{SIR}$) for the $l$th receiver is denoted by $\gamma_l$.

 For $\p=(p_1,\ldots,p_L)\trans\ge \0$, we define the following
 transformation: $\p\mapsto \boldsymbol{\gamma}(\p)$, where
 \begin{equation}\label{defsp}
 \gamma_l(\p):=\frac{g_{ll}p_l}{\sum_{j\ne l}g_{lj}p_j+n_l},
 \;l=1,\ldots,L,
 \; \gam(\p)=(\gamma_1(\p),\ldots,\gamma_L(\p))\trans.
 \end{equation}

 Define
 \begin{equation}
 \label{matrixF} F=[f_{ij}]_{i,j=1}^L, \textrm{ where }
 f_{ij}=\left\{\begin{array}{cl}
    0, & \mbox{if} \,\, i = j \\
    \frac{g_{ij}}{g_{ii}}, & \mbox{if} \,\, i \ne j
       \end{array}\right.
 \end{equation}
 and
 \begin{equation}\label{defbv}
 \g=(g_{11},\ldots,g_{LL})\trans,\;\n=(n_1,\ldots,n_L)\trans,
 \mathbf{v}= \displaystyle \left(\frac{n_1}{g_{11}},
 \frac{n_2}{g_{22}}, \dots,\frac{n_L}{g_{LL}} \right)\trans.
 \end{equation}
 Then
 \begin{equation}\label{ndefgamp}
 \gam(\p)=\p\circ(F\p+\v)^{-1}.
 \end{equation}

 \begin{claim}\label{invtran}  Let $\p\ge \0$ be a
 nonnegative vector.  Assume that $\gam(\p)$ is defined by (\ref{defsp}).
 Then $\rho(\diag(\gam(\p))F)<1$, where $F$ is defined by
 (\ref{matrixF}).  Hence, for $\gam=\gam(\p)$,
 \begin{equation}\label{invgamf}
 \p=P(\gam):=(I-\diag (\gam) F)^{-1}\diag(\gam)\v.
 \end{equation}
 Vice versa, if $\gam$ is in the set
 \begin{equation}\label{gamcond}
 \Gam:=\{\gam\ge \0,\;\rho(\diag (\gam) F) <1\},
 \end{equation}
 then the vector $\p$ defined by
 (\ref{invgamf}) is nonnegative.  Furthermore, $\gam(P(\p))=\gam$.
 That is, $\gam: \R_+^L \to \Gam$, and $P:\Gam\to \R_+^L$ are
 inverse mappings.
 \end{claim}

 \proof
 Observe that (\ref{defsp}) is equivalent to the
 equality
 \begin{equation}\label{pgamid}
 \p=\diag(\gam)F\p+\diag(\gam)\v.
 \end{equation}
 Assume first that $\p$ is a positive vector, i.e.,
 $\p> \0$.  Hence, $\gam(\p)> \0$.
 Since all off-diagonal entries of $\F$ are positive it follows
 that the matrix $\diag(\gam)F$ is irreducible.
 As $\v>\0$, we deduce that $\max_{i\in
 [1,n]}\frac{(\diag(\gam)F\p)_i}{p_i}<1$.  The $\min\max$
 characterization of Wielandt of $\rho(\diag(\gam)F)$,
 \cite{Wie} and
 \cite[(38), pp.64]{Gan}, implies $\rho(\diag(\gam)F)<1$.
 Hence, $\gam(\p)\in\Gam$.  Assume now that $\p\ge\0$.
 Note that $p_i>0\iff \gamma_i(\p)>0$.
 So $\p=\0\iff \gamma(\p)=0$.  Clearly,
 $\rho(\gam(\0)F)=\rho(0_{L\times L})=0<1$.  Assume now that $ \p\gneq\0$.
 Let $\cA=\{i:\; p_i>0\}$.  Denote $\gam(\p)(\cA)$ the vector
 composed of positive entries of $\gam(\p)$.
 Let $F(\cA)$ be the principal
 submatrix of $\F$ with rows and columns in $\cA$.
 It is straightforward to see that
 $\rho(\diag(\gam(\p))F)=\rho(\diag(\gam(p)(\cA)F(\cA))$.
 The arguments above imply that
 $$\rho(\diag(\gam(\p))F)=\rho(\diag(\gam(\p)(\cA)\F(\cA))<1.$$

 Assume now that $\gam\in\Gam$.  Then
 \begin{equation}\label{neumexp}
 (I-\diag(\gam)F)^{-1}=\sum_{k=0}^{\infty}
 (\diag(\gam)F)^k\ge 0_{L\times L}.
 \end{equation}
 Hence, $P(\gam)\ge \0$.  The definition of $P(\gam)$ implies
 that $\gam(P(\gam))=\gam$.
 \qed

 \begin{claim}\label{monotcl}  The set $\Gam\subset \R_+^L$ is
 monotonic with respect to the order $\ge$.  That is if $\gam\in\Gam$
 and $\gam \ge\bet \ge \0$ then $\bet \in \Gam$.
 Furthermore, the function $P(\gam)$ is monotone on $\Gam$.
 \begin{equation}\label{monotP}
 P(\gam)\ge P(\bet) \textrm{ if } \gam\in\Gam \textrm{ and }
 \gam \ge \bet \ge \0.
 \end{equation}
 Equality holds if and only if $\gam=\bet$.
 \end{claim}
 \proof  Clearly, if $\gam \ge \bet \ge\0$ then
 $\diag(\gam)F\ge \diag(\bet)F$ which implies
 $\rho(\diag(\gam)F)\ge \rho(\diag(\bet)F)$.
 Hence, $\Gam$ is monotonic.
 Use the Neumann expansion (\ref{neumexp}) to deduce the
 monotonicity of $P$.  The equality case is straightforward.
 \qed

 Note that $\gam(\p)$ is not monotonic in $\p$.
 Indeed, if one increases only the $i$th coordinate of $\p$,
 then one increases the $i$th coordinate of $\gam(\p)$ and
 decreases all other coordinates of $\gam(\p)$.

 As usual, let $\e_i=(\delta_{i1},\ldots,\delta_{iL})\trans, \;i=1,\ldots,L$ be
 the standard basis in $\R^L$.  In what follows, we need the
 following result.

 \begin{theo}\label{descgampset} Let $l\in[1,L]$ be an
 integer and $a>0$.  Denote $[0,a]_l\times \R_+^{L-1}$ the set
 of all $\p=(p_1,\ldots,p_L)\trans\in\R_+^L$ satisfying $p_l\le
 a$.  Then the image of the set $[0,a]_l\times \R_+^{L-1}$
 by the map $\gam$ (\ref{defsp}), is given by
 \begin{equation}\label{descgampsetl}
 \rho(\diag(\gam)(F+ (1/a)\v\e_l\trans))\le 1,\; \0 \le\gam.
 \end{equation}
 Furthermore, $\p=(p_1,\ldots,p_L)\in\R_+^L$ satisfies the
 condition $p_l=a$ if and only if $\gam=\gam(\p)$ satisfies
 \begin{equation}\label{eqpacase}
 \rho(\diag(\gam)(F+(1/a)\v\e_l\trans))=1.
 \end{equation}
 \end{theo}
 \proof  Suppose that $\gam$ satisfies (\ref{descgampsetl}).
 We claim that $\gam\in\Gam$.
 Suppose first that
 $\gam>\0$.  Then $\diag(\gam)(F+t_1\v\e_l\trans)\lneq \diag(\gam)(F+t_2\v\e_l\trans)$
 for any $t_1 <t_2$.
 \cite[Lemma 2, \S2, Ch. XIII]{Gan} yields
 \begin{eqnarray}\label{incrspect}
 \rho(\diag(\gam)F)< \rho(\diag(\gam)(F+t_1\v\e_l\trans))<
 \rho(\diag(\gam)(F+t_2\v\e_l\trans))<\\\rho(\diag(\gam)(F+(1/a)\v\e_l\trans))
 \le 1 \textrm{ for } 0<t_1<t_2< 1/a. \nonumber
 \end{eqnarray}
 Thus $\gam\in \Gam$.
 Combine the above argument with the arguments
 of the proof of Claim \ref{invtran} to deduce that $\gam\in
 \Gam$ for $\gam\ge \0$.

 We now show that $P(\gam)_l\le a$.  The continuity of $P$
 implies that it is enough to consider the case $\gam>\0$.
 Combine the Perron-Frobenius theorem with (\ref{incrspect}) to
 deduce
 \begin{equation}\label{fundinl}
 0 < \det(I-\diag(\gam)(F+t\v\e_l\trans)) \textrm{ for } t\in [0,a^{-1}).
 \end{equation}
 We now expand the right-hand side of the above inequality.
 Let $B=\x\y\trans\in \R^{L\times L}$ be a rank one matrix.
 Then $B$ has $L-1$ zero eigenvalues and one eigenvalue equal
 to $\y\trans \x$.  Hence, $I-\x\y\trans$ has
 $L-1$ eigenvalues equal to $1$ and one eigenvalue is $(1-\y\trans\x)$.
 Therefore, $\det(I-\x\y\trans)=1-\y\trans \x$.
 Since $\gam\in\Gam$ we get that $(I-\diag(\gam)F)$ is
 invertible.  Thus, for any $t\in\R$
 \begin{eqnarray}\nonumber
 \det(I-\diag(\gam)(F+t\v\e_l\trans))=\\
 \det(I-\diag(\gam)F)\det(I-t((I-\diag(\gam)F)^{-1}
 \diag(\gam)\v)\e_l\trans)\label{detexpnl}\\
 \det(I-\diag(\gam)F)(1 -t\e_l\trans
 (I-\diag(\gam)F)^{-1}\diag(\gam)\v)\nonumber.
 \end{eqnarray}
 Combine (\ref{fundinl}) with the above identity
 to deduce
 \begin{equation}\label{plineq}
 1> t\e_l\trans
 (I-\diag(\gam)F)^{-1}\diag(\gam)\v=tP(\gam)_l \textrm{ for }
 t\in [0,a^{-1}).
 \end{equation}
 Letting $t\nearrow a^{-1}$, we deduce that $P(\gam)_l\le a$.
 Hence, the set of $\gam$ defined by (\ref{descgampsetl}) is a
 subset of $\gam([0,a]_l\times \R_+^{L-1})$.

 Let $\p\in [0,a]_l\times \R_+^{L-1}$ and denote
 $\gam=\gam(\p)$.  We show that $\gam$ satisfies
 (\ref{descgampsetl}).  Claim \ref{invtran} implies that
 $\rho(\diag(\gam)F)<1$.  Since $\p=P(\gam)$ and $p_l\le a$ we
 deduce (\ref{plineq}).  Use (\ref{detexpnl}) to deduce
 (\ref{fundinl}).  As $\rho(\diag(\gam)F)<1$, the inequality (\ref{fundinl})
 implies that $\rho(\diag(\gam)F+t\v\trans\e_l)<1$ for $t\in
 (0, a^{-1})$.  Hence, (\ref{descgampsetl}) holds.

 It is left to show the condition (\ref{eqpacase}) holds if and
 only if $P(\gam)_l=a$.  Assume that
 $\p=(p_1,\ldots,p_L)\trans\in\R^L_+,\;p_l=a$ and let $\gam=\gam(\p)$.
 We claim that equality holds in (\ref{descgampsetl}).
 Assume to the contrary that $\rho(\diag(\gam)(F+(1/a)\v\e_l\trans))< 1$.
 Then, there exists $\bet >\gam$ such that $\rho(\diag(\bet)(F+(1/a)\v\e_l\trans))< 1$.
 Since $P$ is monotonic $P(\bet)_l>p_l=a$.  On the other hand,
 since $\bet$ satisfies (\ref{descgampsetl}), we deduce that
 $P(\bet)_l\le a$.  This contradiction yields (\ref{eqpacase}).  Similarly, if
 $\gam\ge \0$ and (\ref{eqpacase}) then $P(\gam)_l=a$.
 \qed

 \begin{corol}\label{imPbdd}  Let $\bar\p=(\bar p_1,\ldots,\bar p_L)\trans > \0$ be a given
 positive vector.
 Then $\gam([\0,\bar\p])$, the image of the set
 $[\0,\bar\p]$ by the map $\gam$ (\ref{defsp}), is given by
 \begin{equation}\label{descgampset1}
 \rho\left(\diag(\gam)\left(F+(1/\bar\p_l)\v\e_l\trans\right)\right)\le 1, \textrm{ for }
 l=1,\ldots,L, \textrm{ and } \gamma\in\R_+^L .
 \end{equation}
 In particular, any $\gam\in \R_+^L$ satisfying the conditions
 (\ref{descgampset1}) satisfies the inequalities
 \begin{equation}\label{upboundsgam}
 \gam\le \bar\gam=(\bar\gamma_1,\ldots,\bar\gamma_L)\trans, \textrm{ where }
 \bar\gamma_l =\frac{\bar p_l}{v_l},\; i=1,\ldots,L .
 \end{equation}

 \end{corol}
 \proof  Theorem \ref{descgampset} yields that  $\gam([\0,\bar\p])$
 is given by (\ref{descgampset1}).
 (\ref{ndefgamp}) yields
 $$\gamma_l(\p)=\frac{p_l}{((F\p)_l+v_l)}\le \frac{p_l}{v_l}\le
 \frac{\bar p_l}{v_l} \textrm{ for } \p\in [\0,\bar\p].$$
 Note that equality holds for $\p=\bar p_l\e_l$.  \qed

 \section{The sum rate maximization problem}
 \label{sumratesection}
 We assume a singe-user decoder at each user, i.e.,
 treating interference as additive Gaussian noise, and all users has perfect channel state information at the receiver.
We also assume that fading occurs sufficiently slowly in the channel, i.e., flat-fading, so that the channel can be considered essentially fixed during transmission.  Assuming that all users employ achievable random Gaussian codes, we can use the Shannon capacity formula for the maximum information rate on a link. In practice, Gaussian codes can be replaced by finite-order signal constellations such as the use of quadrature-amplitude modulation (QAM) or other practical (suboptimal) coding schemes. Assuming a fixed bit error rate (BER) at the receiver, the Shannon capacity formula can be used to deduce the achievable data rate of the $l$th user as \cite{Cover91}:
\begin{equation}
\label{thruput}
\log \left(1+\frac{\gamma_{l}(\mathbf{p})}{\Gamma} \right) \quad \mbox{nats/symbol},
\end{equation}
where $\Gamma$ is the $\mathsf{SNR}$ gap to capacity, which is always greater than 1. In this paper, we absorb $(1/\Gamma)$ into $G_{ll}$ for all $l$, and write the achievable data rate as $\log (1+\gamma_{l}(\mathbf{p}))$.

Let $\w=(w_1,\ldots,w_L)\trans\ge 0$ be a
 given probability vector, which is assigned by the network to each link (to reflect some long-term priority).
 The problem of maximizing sum
 rate in a Gaussian interference-limited channel can be stated as the
 following maximum problem
 \begin{equation}
 \label{maxprob}
 \max_{\0\le\p\le\bar\p} \sum_{l=1}^L w_l \log(1+\gamma_l(\p)).
 \end{equation}
 \begin{lemma}\label{fnescond}  Let $\w$ be a probability vector,
 and assume that $\p^{\star}=(p_1^{\star},\ldots,p_L^{\star})\trans$
 is a maximal solution to (\ref{maxprob}).  Then
 $p_i^{\star}=\bar p_i$ for some $i$.  Furthermore if $w_j=0$
 then $p_j^{\star}=0$.
 \end{lemma}
 \proof  Assume to the contrary that $\p^{\star} < \bar\p$.
 Let $\gam^{\star}=\gam(\p^{\star})$.  Since $P$ is continuous on
 $\Gam$, there exists $\gam\in\Gam$ such that $\gam >
 \gam^{\star}$ such that $P(\gam) < \bar\p$.  Clearly,
 $\Phi_{\w}(\gam(\p^{\star}))< \Phi_{\w}(\gam)$.
 As $\gam=\gam(P(\gam))$, we deduce that $\p^{\star}$ is not a
 maximal solution to (\ref{maxprob}), contrary to our
 assumptions.

 Suppose that $w_j=0$.  For $\p=(p_1,\ldots,p_L)\trans$, let
 $\p_j$ be obtained from $\p$ by replacing the $jth$ coordinate
 in $\p$ by $0$.  Assume that $p_j>0$.   Then
 $\gamma_i(\p)<\gamma_i(\p_j)$ for $i\ne j$. Since $w_j=0$,
 it follows that
 $\Phi_{\w}(\gam(\p))< \Phi_{\w}(\gam(\p_j))$.
 \qed

 Combine the above lemma with Theorem \ref{descgampset} and
 Corollary \ref{imPbdd} to deduce an alternative formulation of
 (\ref{maxprob}).
 \begin{theo}\label{altformp}  The maximum problem
 (\ref{maxprob}) is equivalent to the maximum problem.
 \begin{equation}
 \label{nonconvex1}
 \begin{array}
 [c]{rl}
 \mbox{maximize} & \sum_l w_l \log(1+\gamma_l)\\
 \mbox{subject to} &  \rho(\diag(\boldsymbol{\gamma})
 (F+(1/\bar{p}_l)\mathbf{v}\mathbf{e}\trans_l)) \le 1 \,\;\; \forall \, l\in\an{L}, \\
 \mbox{variables:} & \gamma_l, \,\;\; \forall \, l.
 \end{array}
 \end{equation}
 $\gam^{\star}$ is a maximal solution of the above problem if
 and only if $P(\gam^{\star})$ is a maximal solution
 $\p^{\star}$ of the problem (\ref{maxprob}).  In particular,
 any maximal solution $\gam^{\star}$ satisfies the equality
 (\ref{descgampset1}) for some integer $l\in [1,L]$.
 \end{theo}

 We now give the following simple necessary conditions for a
 maximal solution $\p^{\star}$ of (\ref{maxprob}).
 We first need the
 following result, which is obtained by straightforward
 differentiation.
 \begin{lemma}\label{gradhas}
 Denote by
 \[\nabla
 \Phi_{\w}(\gam)=\left(\frac{w_1}{1+\gamma_1},\ldots,\frac{w_L}{1+\gamma_L}\right)\trans
 =\w\circ (\1+\gam)^{-1}\]
 the gradient of $\Phi_{\w}$.  Let $\gam(\p)$ be defined as in
 (\ref{defsp}).  Then $\rH(\p)=[\frac{\partial \gamma_i}{\partial p_j}]
 _{i=j=1}^L $,  the Hessian matrix of $\gam(\p)$, is given by
 \[H(\p)=\diag((F\p+\v)^{-1})(-\diag(\gam(\p))F+I).\]
 In particular,
 \[\nabla_{\p} \Phi_{\w}(\gam(\p))=\rH(\p)\trans \nabla
 \Phi_{\w}(\gam(\p)).\]
 \end{lemma}

 \begin{corol}\label{neccondmaxsol}
 Let $\p^{\star}=(p_1^{\star},\ldots,p_L^{\star})\trans$ be a maximal solution to
 the problem (\ref{maxprob}).  Divide the set
 $\an{L}=\{1,\ldots,L\}$ to the following three disjoint sets
 $\rS_{\max},\rS_{\inn},\rS_0$:
 \[
 \rS_{\max}=\{i\in \an{L},\;p_i^{\star}=\bar p_i\},
 \;\rS_{\inn}=\{i\in\an{L},\;p_i^{\star}\in (0,\bar p_i)\},
 \;\rS_0=\{i\in\an{L},\; p_i^{\star}=0\}.
 \]
 Then the following conditions hold.
 \begin{eqnarray}
 &&(\rH(\p^{\star})\trans \nabla \Phi_{\w}(\gam(\p^{\star})))_i\ge 0 \textrm{ for
 } i\in \rS_{\max},\nonumber\\
 &&(\rH(\p^{\star})\trans \nabla \Phi_{\w}(\gam(\p^{\star})))_i= 0 \textrm{ for
 } i\in \rS_{\inn},\label{neccondmaxsol1}\\
 &&(\rH(\p^{\star})\trans \nabla \Phi_{\w}(\gam(\p^{\star})))_i\le 0 \textrm{ for
 } i\in \rS_{0}\nonumber
 \end{eqnarray}

 \end{corol}
 \proof  Assume that $p_i^{\star}=\bar p_i$.  Then
 $\frac{\partial}{\partial p_i} \Phi_{\w}(\gam(\p))(\p^{\star})\ge 0$.
 Assume that $0<p_i^{\star}<\bar p_i$.  Then
 $\frac{\partial}{\partial p_i} \Phi_{\w}(\gam(\p))(\p^{\star})=0$.
 Assume that $p_i^{\star}=0$.  Then
 $\frac{\partial}{\partial p_i} \Phi_{\w}(\gam(\p))(\p^{\star})\le 0$.
 \qed

 We now show that the maximum problem (\ref{nonconvex1}) can be
 restated as the maximum problem of convex function on a closed
 unbounded domain.
 For $\gam=(\gamma_1,\ldots,\gamma_L)\trans >0$ let $\tilde\gam=\log\gam$, i.e.
 $\gam=e^{\tilde\gam}$.
 Recall that for a nonnegative irreducible matrix
 $B\in\R_+^{L\times L}$ $\log\rho(e^{\x} B)$ is a convex function \cite{Kin}.
 Furthermore, $\log (1+e^t)$ is a strict convex function in $t\in \R$.
 Hence, the maximum problem (\ref{nonconvex1}) is equivalent to the
 problem

 \begin{equation}
 \label{newmax2}
 \begin{array}
 [c]{rl}
 \mbox{maximize} & \sum_l w_l \log(1+e^{\tilde{\gamma}_l})\\
 \mbox{subject to} & \log \rho(\diag(e^{\tilde{\gam}})
 (F+(1/\bar{p}_l)\mathbf{v}\mathbf{e}\trans_l)) \le 0 \,\;\; \forall \, l\in\an{L}, \\
 \mbox{variables:} & \tilde{\gam}=(\tilde{\gamma}_1,\ldots,\tilde{\gamma}_n)\trans \in \R^L.
 \end{array}
 \end{equation}
 The unboundedness of the convex set in (\ref{newmax2})
 is due to the identity $0=e^{-\infty}$.
 In view of Lemma \ref{fnescond}, it is enough to consider the
 maximal problem (\ref{maxprob}) in the
 case where $\w>\0$.
 \begin{theo}\label{locmax}  Let $\w>\0$ be a probability vector.
 Consider the maximum problem (\ref{maxprob}).  Then any
 point $\0\le \p^{\star}\le \bar\p$ satisfying the conditions
 (\ref{neccondmaxsol1}) is a local maximum.
 \end{theo}
 \proof
 Since $\w>\0$, $\Phi_{\w}(e^{\tilde\gam})$ is a
 strict convex function in $\tilde\gam\in\R^L$.
 Hence, the maximum of (\ref{newmax2}) is achieved exactly on
 the extreme points of the closed  unbounded set specified in
 (\ref{newmax2}).  (It may happen that some coordinate of the
 extreme point are $-\infty$.)  Translating this observation to
 the maximal problem (\ref{maxprob}), we deduce the theorem.
 \qed

  We now give simple lower and upper bounds on the value of
 (\ref{maxprob}).

 \begin{lemma}\label{lowupbound}  Consider the maximal problem
 (\ref{maxprob}).  Let
 $B_l=(F+(1/\bar{p}_l)\mathbf{v}\mathbf{e}\trans_l))$  for
 $l=1,\ldots,L$.  Denote $R=\max_{l\in\an{L}} \rho(B_l)$.
 Let $\bar\gam$ be defined by (\ref{upboundsgam}).  Then
 \[
 \Phi_{\w}((1/R)\1)\le\max_{\p\in [\0,\bar\p]}\Phi_{\w}(\gam(\p)\le
 \Phi_{\w}(\bar\gam).
 \]

 \end{lemma}
 \proof  By Corollary \ref{imPbdd}, $\gam(\p)\le \bar\gam$ for
 $\p\in [\0,\bar\p]$.  Hence, the upper bounds holds.
 Clearly, for $\gam=(1/R)\1$, we have that
 $\rho(\diag(\gam)B_l)\le 1$ for $l\in \an{L}$.
 Then, from Theorem \ref{altformp}, $\Phi_{\w}((1/R)\1)$ yields the lower bound. Equality is achieved in the lower bound when $\mathbf{p}^{\star} = t \mathbf{x}(B_{i})$, where $i=\arg \max_{l\in\an{L}} \rho(B_l)$, for some $t>0$.
 \qed

 \section{Relaxations of the maximal problem}

In this section, we study several relaxed versions of (\ref{nonconvex1}), which will be used later to construct algorithms to solve (\ref{maxprob}).

 From the definition of $\boldsymbol{\gamma}(p)$, we deduce that
 \begin{equation}\label{peq1}
 p_l=\gamma_l(\p) \left(\frac{n_l}{g_{ll}}+\sum_{j\ne
 l}\frac{g_{lj}}{g_{ll}}p_j \right).
 \end{equation}
 Define
 \begin{equation}\label{defF}
 \tilde{F}=[\tilde{f}_{lj}]_{l,j=1}^L,\; \tilde{f}_{ll}=\frac{n_l}{g_{ll}\bar p_l},
 \;\tilde{f}_{lj}=\frac{g_{lj}}{g_{ll}} \textrm{ for } j\ne l, \;
 l=1,\ldots,L.
 \end{equation}

 \begin{lemma}\label{peqlem}  Let $\0\le\p\le\bar\p$.  Assume
 that $\boldsymbol{\gamma}(\p)$ and $\tilde{F}$ are defined by (\ref{defsp}) and
 (\ref{defF}), respectively.  Then
 \begin{equation}\label{peq2}
 \p\ge \diag(\boldsymbol{\gamma}(\p)) \tilde{F}\p,
 \end{equation}
 and
 \begin{equation}\label{peq3}
 \rho(\diag (\boldsymbol{\gamma}(\p))\tilde{F})\le 1.
 \end{equation}
 \end{lemma}
 \proof
 The assumption that $0\le p_l\le \bar p_l$ implies that
 $\frac{n_l}{g_{ll}}\ge \tilde{f}_{ll} p_l$.
 Then using (\ref{peq1}), the definition of $\tilde{F}$ and the above observation implies
 (\ref{peq2}).  The inequality (\ref{peq3}) is a consequence
 of Wielandt's characterization of the spectral radius of an
 irreducible matrix \cite{Gan}.  Indeed, if $\p>\0$, i.e. all
 the coordinates of $\p$ are positive, then $\boldsymbol{\gamma}(\p)>0$.  Hence,
 $\diag(\boldsymbol{\gamma}(\p)))\tilde{F}$ is a positive matrix.  Then Wielandt's
 characterizations claims
 $$\rho(\diag(\boldsymbol{\gamma}(\p))\tilde{F})\le \max_{l=1,\ldots,L}
 \frac{(\diag(\boldsymbol{\gamma}(\p))\tilde{F}\p)_l}{p_l}\le 1.$$
 Observe next that if $p_l=0$, then $\boldsymbol{\gamma}(\p)_l=0$.
 So if some of $p_l=0$, then $\rho(\diag(\boldsymbol{\gamma}(\p))\tilde{F})$ is the
 spectral radius of the maximal positive submatrix of
 $\diag(\boldsymbol{\gamma}(\p))\tilde{F}$.  Apply to this positive submatrix Wielandt's
 characterization to deduce (\ref{peq3}).  \qed

 \begin{lemma}\label{newmax}  The maximum
 \begin{equation}\label{newmax1}
 \max_{\gam\in\R^L_+, \rho(\diag(\gam)\tilde{F})\le 1}
 \Phi_{\w}(\gam)
 \end{equation}
 is not less than the maximum of the problem (\ref{maxprob}).
 \end{lemma}
 \proof  In view of (\ref{peq3}), we see that the maximum in
 (\ref{newmax1}) is on a bigger set than the maximum in
 (\ref{maxprob}).  \qed

 By considering a change of variable in the constraint set of
 (\ref{newmax1}), we define the following set:
 \begin{equation}\label{defconvset}
 D(\tilde{F})=\{\tilde{\gam}\in\R^{L}, \quad \log\rho(\diag(e^{\tilde{\gam}})\tilde{F})\le
 0\},
 \end{equation}
 which is convex in $\R^L$.

 \begin{corol}
 \label{constrainttight}
 We have $\rho(\diag(\gam^{\prime})\tilde{F}) = 1$ in (\ref{newmax1}), where $\gam^{\prime}$
 solves (\ref{newmax1}) optimally.
 \end{corol}

 \begin{lemma}
 \label{peqlemma}
 If $\p^{\star}=\bar{\mathbf{p}}$ or $\p^{\star}$ is such that
 $p^{\star}_l=0$ for some $l$ and $p^{\star}_j=\bar p_j$ for $j \ne
 l$, then
 \begin{equation}
 \label{peqlemmaeqn}
  \rho(\diag (\boldsymbol{\gamma}(\p^{\star}))\tilde{F}) = 1.
 \end{equation}
 \end{lemma}
 \proof The definition of $\tilde F$ implies (\ref{peqlemmaeqn}) for
 $\p^{\star}=\bar{\mathbf{p}}$.
 Assume now that $p^{\star}_l=0$ for some $l$. Then
 $\gamma_l(\p^{\star})=0$ for some $l$. Then,
 the $lth$ row of $\diag
 (\boldsymbol{\gamma}(\p^{\star}))\tilde{F}$ is zero.
 Let $F_1$ be the submatrix of $F$ obtained by deleting
 the $lth$ row and column.  Let $\gam_1$ be the vector
 obtained from $\gam$ by deleting the $lth$ coordinate.
 Hence, the characteristic polynomial of $\diag(\gam)F$,
 $\det(xI-\diag(\gam)F$, is equal to
 $x\det(xI-\diag(\gam_1)F_1)$.  Therefore, $\rho(\diag(\gam)F)=
 \rho(\diag(\gam_1)F_1)$.  Continuing in this manner, we deduce the
 lemma. \qed

 Corollary \ref{constrainttight} and Lemma \ref{peqlemma} imply that if the optimizer of (\ref{newmax1}) $\gam^{\prime}$ satisfies $P(\gam^{\prime}) \preceq \bar{\mathbf{p}}$, then $P(\gam^{\prime})$ is also the global optimizer of (\ref{maxprob}).

 We now state a stronger relaxation problem than (\ref{newmax1}).
 Consider the transformation $\gam(\p,\n)$ given in (\ref{defsp}) as
 a function of $\p$ and the noise $\n>\0$.  Clearly,
 $\gam(\p,\n_1)\gneq\gam(\p,\n_2)$ for any $\n_2\gneq \n_1>\0$ and
 $\p\ge \0$.  Hence, it would be useful to consider the limiting
 case $\delta(\p):=\gam(\p,\0)$:
 \begin{equation}\label{defbetp}
 \beta_l(\p):=\frac{g_{ll}p_l}{\sum_{j\ne l}g_{lj}p_j},
 \;l=1,\ldots,L,
 \; \bet(\p)=(\beta_1(\p),\ldots,\beta_L(\p))\trans, \p \gneq \0.
 \end{equation}
 Hence, $\beta(\p)$ can be viewed as communication in \emph{noiseless
 channels} \cite{Aein73}.   Another way to look at the noiseless channels is to allow
 channels with unlimited power $\p$, i.e. $\bar p_i=\infty,
 i=1,\ldots,L$.  This remark is implied by the identity
 $$\gam(t\p,\n)=\gam(\p, \frac{1}{t}\n) \textrm{ for any } t>0.$$
 Denote by $\Pi_L=\{\p\in \R^L_+,\;\1\trans \p=1\}$ the set of
 probability vectors in $\R_+^L$.
 Since $\beta(\p)=\beta(t\p)$ for any $t>0$, we deduce
 \begin{equation}\label{maxclchan}
 \sup_{\p\in \R^L_{+} \backslash\{\0\}}
 \Phi_{\w}(\bet(\p))=\max_{\p\in\Pi_L} \Phi_{\w}(\bet(\p)).
 \end{equation}
 Clearly, the maximum in (\ref{maxclchan}) is greater than the
 maximum in (\ref{maxprob}).

 Note that the definition of $F$ given by (\ref{matrixF}) implies
 that
 \begin{equation}\label{betFid}
 \diag(\bet(\p))F\p=\p, \textrm{ and } \rho(\diag(\bet(\p))F)=1 \textrm{ for
 any } \p\gneq \0.
 \end{equation}
 Hence, the maximum (\ref{maxclchan}) is equal to
 \begin{equation}\label{maxclchan1}
 \max_{\bet\in \R_+^L, \rho(\diag(\bet) F)\le 1} \Phi_{\w}(\bet)=
 \max_{\bet\in \R_+^L, \rho(\diag(\bet) F)= 1} \Phi_{\w}(\bet).
 \end{equation}
 Since $F\lneq \tilde F$ and $\tilde F$ is positive, it follows that
 $\rho(\diag(\gam) F)< \rho(\diag(\gam) \tilde F)$ for any $\gam \gneq \0$ \cite{Gan}.
 Hence, the maximum in
 (\ref{maxclchan1}) is greater than the maximum in (\ref{newmax1}).

Since it is generally difficult to determine precisely the spectral radius of a given matrix \cite{Varga63}, the relaxed problems given by (\ref{newmax1}) and (\ref{maxclchan1}) enables one to find an upper bound to (\ref{maxprob}) quickly when $L$ is fairly large.

 We next modify the function $\Phi_{\w}(e^{\tilde\gam})$
 appearing in (\ref{maxapproxD1}).  Consider the function
 $\log(1+e^x)$.  If $e^x$ is relatively small, then $\log (1+e^x)$ is close
 to zero.  If $e^x$ is assumed to be bigger than 1, then $\log(1+e^x)\approx \log e^x=x$.
 Hence, it is reasonable to approximate $\log(1+e^x)$ by $x$ (Note that $\log(1+e^x)>x$).
 \begin{theo}\label{tilfrelsimp}  Let $\w=(w_1,\ldots,w_L)\trans$ be a
 positive probability vector.  Then
 \begin{equation}\label{tilfrelsimp1}
 \max_{\gam=(\gamma_1,\ldots,\gamma_L)\trans>\0,
 \rho(\diag(\gam) \tilde F)\le 1} \sum_{l=1}^L w_l
 \log\gamma_l= \sum_{l=1}^L w_l \log \gamma_l^{\star},
 \end{equation}
 where
 $\gam^{\star}=(\gamma_1^{\star},\ldots,\gamma_L^{\star})\trans
 >\0$ is the unique vector satisfying the following conditions.
 Let $A^{\star}=\diag(\gam^{\star})\tilde F$.  Then
 $\rho(A^{\star})=1$ and $\x(A^{\star})\circ \y(A^{\star})=\w$.

 \end{theo}
 \proof  In the maximum problem (\ref{tilfrelsimp1}), it is
 enough to consider $\gam>\0$ such that $ \rho(\diag(\gam) \tilde
 F)=1$.  Then each of such $\gam$ is of the form $\p\circ
 (F\p)^{-1}$, for some $\0<\z\;(=\x(\diag(\gam) \tilde F)$ in Theorem \ref{apfkineq}.
 Use Theorem \ref{apfkineq} and Corollary \ref{corthm3.2fk} to
 deduce the theorem.  \qed

 Combining Theorem \ref{apfkineq} and Corollary
 \ref{corthm3.2ofdiag} we deduce the following result.
 \begin{theo}\label{frelsimp}  Let $\w=(w_1,\ldots,w_L)\trans$ be a
 positive probability vector satisfying the conditions (\ref{majcon}).
 Let $\hat F$ be either equal to $F$ or to $F+(1/\bar
 p_l)\v\e_l\trans$ for some $l\in\an{L}$, respectively.
 Then
 \begin{equation}\label{frelsimp1}
 \max_{\gam=(\gamma_1,\ldots,\gamma_L)\trans>\0,
 \rho(\diag(\gam) \hat F)\le 1} \sum_{l=1}^L w_l
 \log\gamma_l= \sum_{l=1}^L w_l \log \gamma_l^{\star},
 \end{equation}
 where
 $\gam^{\star}=(\gamma_1^{\star},\ldots,\gamma_L^{\star})\trans
 >\0$ is a vector satisfying the following conditions.
 Let $A^{\star}=\diag(\gam^{\star}) \hat F$.  Then
 $\rho(A^{\star})=1$ and $\x(A^{\star})\circ \y(A^{\star})=\w$.

 \end{theo}
 Note that the choice $\hat F=F+(1/\bar p_l)\v\e_l\trans$
 applies, if we assume that the optimal solution
 $\p^{\star}=(p_1^{\star},\ldots,p_L^{\star})\trans$ of
 the  maximal problem (\ref{maxprob}) satisfies
 $p_l^{\star}=\bar p_l$ and $p_j^{\star}<\bar p_j$ for $j\ne
 l$.

 The above last two theorems enable us to choose $\w$ for
 which  we know the solution to the maximal problems
 (\ref{tilfrelsimp1}) and (\ref{frelsimp1}).  Namely,
 choose $\bet_1, \bet_2 \ >\0$ such that $A_1=\diag(\bet_1)\tilde
 F,A_2=\diag(\bet_2)\hat F$ have spectral radius one.  Let
 $\w_i=\x(A_i)\circ \y(A_i)$ for $i=1,2$.
 Then for $\w_1$, (\ref{tilfrelsimp1}) has the unique maximal
 solution $\gam^{\star}=\bet_1$.  For $\w_2$, (\ref{frelsimp1}) has a maximal
 solution $\gam^{\star}=\bet_2$.  Note that in view of Theorem \ref{apfkineq},
 $\w_2$ does not have to satisfy the conditions (\ref{majcon}).

 \section{Algorithms for sum rate maximization}
\label{algosection}
 In this section, we outline three algorithms for finding and
 estimating the maximal sum rates.
 In view of Lemma \ref{fnescond}, it is enough to consider the
 maximal problem (\ref{maxprob}) in the
 case where $\w>\0$. Theorem \ref{locmax}
 gives rise to the following algorithm,
 which is the gradient algorithm in the variable $\p$
 in the compact polyhedron $[\0,\bar\p]$.
 \begin{algo}\label{gradalgo}
 $ $
 \begin{enumerate}
 \item  Choose $\p_0\in [\0,\bar\p]$:
 \begin{enumerate}
 \item Either at random;
 \item or $\p_0=\bar\p$.
 \end{enumerate}
 \item Given $\p_k=(p_{1,k},\ldots,p_{L,k})\trans\in [\0,\bar\p]$ for $k\ge 0$, compute
 $\a=(a_1,\ldots,a_L)\trans=\nabla_{\p} \Phi_{\w}(\gam(\p_k))$.  If $\a$ satisfies the
 conditions (\ref{neccondmaxsol1}) for $\p^{\star}=\p_k$, then $\p_k$ is the output.
 Otherwise let $\b=(b_1,\ldots,b_L)\trans$ be defined as
 follows.
 \begin{enumerate}
 \item $b_i=0$ if $p_{i,k}=0$ and $a_i<0$;
 \item $b_i=0$ if
 $p_{i,k}=\bar p_i$ and $a_i>0$;
 \item
 $b_i=a_i$ if $0<p_i<\bar p_i$.
 \end{enumerate}
 Then $\p_{k+1}=\p_k+t_k \b$, where $t_k>0$ satisfies the conditions
 $\p_{k+1}\in [\0,\bar\p]$ and $\Phi_{\w}(\gam(\p_k+t_k\b_k))$ increases
 on the interval $[0,t_k]$.

 \end{enumerate}

 \end{algo}

 The problem with the gradient method, and its variations as a conjugate gradient method is that
 it is hard to choose the optimal value of $t_k$ in each step, e.g. \cite{Avr}.
 We now use the reformulation of the maximal problem
 given by (\ref{newmax2}). Since $\w>\0$, the function
 $\Phi_{\w}(e^{\tilde\gam})$ is
 strictly convex.  Thus, the maximum is achieved only on the
 boundary of the convex set
 \begin{equation}\label{orgconvset}
 D(\{F\})=\{\tilde{\gam} \in\R^{L},\quad \log\rho(\diag(e^{\tilde{\gam}})
 (F+(1/\bar{p}_l)\mathbf{v}\mathbf{e}\trans_l))\le 0, \;\; \forall \; l\}.
 \end{equation}

 If one wants to use numerical methods and
 software for finding the maximum value of convex functions on bounded closed convex sets ,
 e.g., \cite{NoW}, then one needs to consider the
 maximization problem (\ref{newmax2}) with additional constraints:
 \begin{equation}\label{Kcosntr}
 D(\{F\},K)=\{\tilde\gam\in D(\{F\}),\quad \tilde{\gam}\ge -K\1\}.
 \end{equation}
 for a suitable $K\gg 1$.  Note that the above closed set is
 compact and convex.  The following lemma gives the description
 of the set $D(\{F\},K)$.
 \begin{lemma}\label{desDFK}  Let $\bar\p>\0$ be given and let $R$
 be defined as in Lemma \ref{lowupbound}.  Assume that $K>\log R$.  Let
 $\underline{\p}=P(e^{-K}\1)=(e^{K}I-F)^{-1}\v$.  Then
 $D(\{F\},K)\subseteq\log\gam([\underline{\p},\bar\p])$.
 \end{lemma}
  \proof  From the definition of $K$, we have that $e^K>R$.
 Hence, $\rho(e^{-K}B_l)<1$ for $l=1,\ldots,L$.  Thus
 $-K\1\in D(\{F\})$.  Let
 $\underline{\gam}=e^{-K}\1$.  Assume that $\tilde\gam\in
 D(\{F\},K)$.  Then $\tilde\gam\ge -K\1$.  Hence, $\gam=e^{\tilde\gam}\ge
 \underline{\gam}$.  Since $\rho(\diag(\gam)F)<1$, Claim \ref{monotcl}
 yields that $\p=P(\gam)\ge P(\underline{\gam})=\underline{\p}$,
 where $P$ is defined by (\ref{invgamf}).
 The inequality $P(\gam)\le\bar\p$ follows from Corollary
 \ref{imPbdd}.  \qed

 Thus, we can apply the numerical methods to find the maximum of
 the strictly convex function $\Phi_{\w}(e^{\tilde\gam})$ on the closed
 bounded set $D(\{F\},K)$, e.g. \cite{NoW}.
 In particular, we can use the gradient method.  It takes the
 given boundary point $\tilde\gam_k$ to another boundary point of
 $\tilde\gam_{k+1}\in D(\{F\},K)$, in the direction induced by
 the gradient of $\Phi_{\w}(e^{\tilde\gam})$.
 However, the complicated boundary of  $D(\{F\},K)$ will make any
 algorithm expensive.

Furthermore, even though the constraint set in (\ref{nonconvex1}) can be transformed into a strict convex set, it is in general difficult to determine precisely the spectral radius of a given matrix \cite{Varga63}. To make the problem simpler and to enable fast algorithms, we approximate the convex set
 $D(\{F\},K)$ by a bigger polyhedral convex sets as follows.
 Choose a finite number of points $\zet_1,\ldots,\zet_M$
 on the boundary of $D(\{F\})$, which preferably lie in $D(\{F\},K)$.
 Let
 \par\noindent
 $\rH_1(\xit),\ldots,\rH_N(\xit), \xit\in\R^L$ be the $N$
 supporting hyperplanes of $D(\{F\}$.  (Note that we can have
 more than one supporting hyperplane at $\zet_i$, and at most
 $L$ supporting hyperplanes.)  So each $\xit\in D(\{F\},K)$ satisfies the
 inequality $H_j(\xit)\le 0$ for $j=1,\ldots,N$.  Let $\bar\gam$
 be defined by (\ref{upboundsgam}).  Define
 \begin{equation}\label{approxD}
 D(\zet_1,\ldots,\zet_M,K)=\{\xit\in\R^L,\;-K\1\le \xit\le
 \log\bar\gam ,\; \rH_j(\xit)\le 0 \textrm{ for
 } j=1,\ldots, N\}.
 \end{equation}
 Hence, $D(\zet_1,\ldots,\zet_M,K)$ is a polytope which
 contains $D(\{F\},K)$.  Thus
 \begin{eqnarray}\label{maxapproxD1}
 \max_{\tilde\gam\in
 D(\zet_1,\ldots,\zet_M,K)}\Phi_{\w}(e^{\tilde\gam})\ge\\
  \max_{\tilde\gam\in
 D(\{F\},K)}\Phi_{\w}(e^{\tilde\gam}).
 \label{maxapproxD2}
 \end{eqnarray}
 Since $\Phi_{\w}(e^{\tilde\gam})$ is strictly convex, the
 maximum in (\ref{maxapproxD1}) is achieved only at the
 extreme  points of $D(\zet_1,\ldots,\zet_M,K)$.
 The maximal solution can be found using a variant of a
 simplex algorithm \cite{Dan}.  More precisely, one starts at
 some extreme point of $\xit\in D(\zet_1,\ldots,\zet_M,K)$.
 Replace the strictly convex function
 $\Phi_{\w}(e^{\tilde\gam})$ by its first order Taylor expansion $\Psi_{\xit}$ at
 $\xit$.  Then we find another extreme point $\et$ of
 $D(\zet_1,\ldots,\zet_M,K)$, such that
 $\Psi_{\xit}(\et)>\Psi_{\xit}(\xit)=\Phi_{\w}(e^{\xit})$.
 Then we replace $\Phi_{\w}(e^{\tilde\gam})$ by its first order Taylor expansion $\Psi_{\et}$ at
 $\et$ and continue the algorithm.
 Our second proposed algorithm for finding an optimal $\tilde \gam^{\star}$ that maximizes (\ref{maxapproxD1}) is given as follows.
\begin{algo}\label{tayloralgo}
$ $
 \begin{enumerate}
 \item Choose an arbitrarily extreme point $\xit_0 \in D(\zet_1,\ldots,\zet_M,K)$.
 \item \label{tayloralgostep2}
 Let $\Psi_{\xit_k}(\xit) = \Phi_{\w}(e^{\xit_k}) + (\mathbf{w} \circ (\mathbf{1}+e^{\xit_k})^{-1}
 \circ e^{\xit_k})^{\trans}(\xit - \xit_k)$. Solve the linear program $\max_{\xit}
 \Psi_{\xit_k}(\xit)$ subject to $\xit \in D(\zet_1,\ldots,\zet_M,K)$ using the simplex
 algorithm in \cite{Dan} by finding an extreme point $\xit_{k+1}$ of $D(\zet_1,\ldots,\zet_M,K)$,
 such that $\Psi_{\xit_k}(\xit_{k+1})>\Psi_{\xit_k}(\xit_k)=\Phi_{\w}(e^{\xit_k})$.
 \item Compute $\p_{k}=P(e^{\xit_{k+1}})$. If $\p_k \in [0,\bar{\mathbf{p}}]$, compute
 $\a=(a_1,\ldots,a_L)\trans=\nabla_{\p} \Phi_{\w}(\gam(\p_k))$.  If $\a$ satisfies the
 conditions (\ref{neccondmaxsol1}) for $\p^{\star}=\p_k$, then $\p_k$ is the output.
 Otherwise, go to Step \ref{tayloralgostep2} using $\Psi_{\xit_{k+1}}(\xit)$.
 \end{enumerate}
\end{algo}

 As in the previous section, it would be useful to
 consider the following related maximal problem:
  \begin{equation}\label{maxapproxD3}
 \max_{\tilde\gam\in
 D(\zet_1,\ldots,\zet_M,K)}\w\trans\tilde\gam.
 \end{equation}

 This problem given by (\ref{maxapproxD3}) is a standard linear program, which can be solved
 in polynomial time by the classical ellipsoid algorithm
 \cite{Kha}. Our third proposed algorithm for finding an optimal $\tilde \gam^{\star}$ that maximizes (\ref{maxapproxD3}) is given as follows.

\begin{algo}\label{approxalgo}
$ $
 \begin{enumerate}
 \item Solve the linear program $\max_{\tilde\gam} \w\trans\tilde\gam$ subject to $\tilde\gam \in D(\zet_1,\ldots,\zet_M,K)$ using the ellipsoid algorithm in \cite{Kha}.
 \item Compute $\p=P(e^{\tilde\gam})$. If $\p \in [0,\bar{\mathbf{p}}]$, then $\p$ is the output. Otherwise, project $\p$ onto $[0,\bar{\mathbf{p}}]$.
 \end{enumerate}
\end{algo}

We note that $\tilde\gam \in D(\zet_1,\ldots,\zet_M,K)$ in Algorithm \ref{approxalgo} can be replaced by the set of supporting hyperplane $ D(\tilde F,K)=\{\tilde\gam\in \rho(\diag(e^{\tilde\gam})\tilde{F}) \le 1,\quad \tilde{\gam}\ge -K\1\}$ (extending (\ref{defconvset})) or, if $L \ge 3$ and $\mathbf{w}$ satisfies the conditions (\ref{majcon}), $ D( F,K)=\{\tilde\gam\in \rho(\diag(e^{\tilde\gam}) F) \le 1,\quad \tilde{\gam}\ge -K\1\}$ based on the relaxed maximal problems in Section 4. Then Theorem \ref{tilfrelsimp} and \ref{frelsimp} quantify the closed-form solution $\tilde\gam$ computed by Algorithm \ref{approxalgo}.

 We conclude this section by showing how to compute the supporting hyperplanes $\rH_j,
 j=1,\ldots,N$, which define $D(\zet_1,\ldots,\zet_M,K)$.
 To do that, we give a characterization of supporting hyperlanes of
 $D(\{F\})$ at a boundary point $\zet\in\partial D(\{F\})$.
 \begin{theo}\label{descsuphyDF}
 Let $\bar\p=(\bar p_1,\ldots\bar p_L)\trans>\0$ be given.  Consider the convex set
 (\ref{orgconvset}).  Let $\zet$ be a boundary point of $\partial
 D(\{F\})$.  Then $\zet=\log\gam(\p)$, where
 $\0\le \p=(p_1,\ldots,p_L)\trans\le \bar\p$.  The set
 $\cB:=\{l\in\an{L},\; p_l=\bar p_l\}$ is nonempty.  For each
 $B_l=(F+(1/\bar{p}_l)\mathbf{v}\mathbf{e}\trans_l))$ let
 $\rH_l(\xi)$ be defined as in Theorem \ref{suphyp}, where
 $B=B_l$ and $\et=\zet$.  Then $\rH_l\le 0$, for $l\in \cB$, are the
 supporting hyperplanes of $D(\{F\})$ at $\zet$.

 \end{theo}
 \proof  Let $\p=P(e^{\zet})$.
 Theorem \ref{descgampset} implies the set $\cB$ is nonempty.
 Furthermore, $\rho(e^{\zet}B_l)=1$ if and only if
 $p_l=\bar p_l$.  Hence, $\zeta$ lies exactly at the intersection of the
 hypersurfaces $\log\rho(e^{\zet}B_l)=0, l\in \cB$.
 Theorem \ref{suphyp} implies that the supporting hyperplanes
 of $D(\{F\})$ at $\zet$ are $\rH_l(\xit)\le 0$ for $l\in\cB$.
 \qed

 We now show how to choose the boundary points
 $\zet_1,\ldots,\zet_M\in \partial D(\{F\})$ and to compute the
 supporting hyperplanes of $D(\{F\})$ at each $\zeta_i$.
 Let $\underline{\p}=P(e^{-K}\1)=(\underline{p_1},\ldots,\underline{p_L})\trans$
 be defined as in Lemma \ref{desDFK}.  Choose $M_i\ge 2$ equidistant points in
 each interval $[\underline{p}_i,\bar p_i]$.
 \begin{equation}\label{defpji}
 p_{j_i,i}=\frac{j_i\underline{p}_i+(M_i-j_i)\bar p_i}{M_i} \textrm{ for
 } j_i=1,\ldots, M_i, \textrm{ and } i=1,\ldots,L.
 \end{equation}
 Let
 $$\cP=\{\p_{j_1,\ldots,j_L}=(p_{j_1,1},\ldots,p_{j_L,L})\trans,\;\min(\bar
 p_1-p_{j_1,1},\ldots,\bar p_{L}-p_{j_L,L})=0\}.$$
 That is, $\p_{j_1,\ldots,j_L}\in\cP$ if and
 only $\p_{j_1,\ldots,j_L}\nless \bar\p$.
 Then
 $$\{\zet_1,\ldots,\zet_M\}=\log \gam(\cP).$$
 The supporting hyperplanes of $D(\{F\})$ at each $\zet_i$
 are given by Theorem \ref{descsuphyDF}.

\section{Extension to DSL channels}
In this section, we extend our previous results in Section \ref{sumratesection} to the multiuser DSL channel, where a common spectrum is divided into $K$ frequency tones denoted by $\an{K}$. For simplicity, we assume the standard synchronous discrete multi-tone (DMT) modulation over DSL, where orthogonality among subchannels of the intended signal and the subchannels of the interference signal in different frequency tones is maintained. Thus, transmissions can be modeled independently on each tone. The achievable rate at tone $k$ can be modeled as
\begin{equation}
\label{ratedsl}
\log\left( 1+\frac{g_{ll}^{k} p_l^{k}}{\sum_{j \ne l}g_{lj}^{k} p_j^{k}+n_l^{k}} \right).
\end{equation}

The total data rate for each user is then obtained by adding its transmitted bits over all the $K$ tones.
The total power budget of the $l$th user is constrained (across all $K$ tones) by
\begin{equation}
\sum_{k=1}^{K} p_l^{k} \le \bar{p}_l.
\end{equation}
It is easy to see that our previous model in Section \ref{modelsection} is a special case of the DSL model assuming standard synchronous discrete multi-tone (DMT) modulation when $K=1$.

For brevity of notations, we define $\mathbf{p} \in \R^{(K \times L)}$ as a vector that stacks the $K \times L$ power allocation lined-up according to tones of all users. For example, if $L=2$ and $K=2$, then $\mathbf{p}=(p^{1}_{1},p^{2}_{1},p^{1}_{2},p^{2}_{2})^{\trans}$. Similarly, we define $\boldsymbol{\gamma} \in \R^{(K \times L)}$ as the $\mathsf{SIR}$ allocation, the following matrix
 \begin{equation}
 \label{matrixFdsl} F=[f^{k}_{ij}]_{i,j=1}^{K \times L}, \textrm{ where }
 f^{k}_{ij}=\left\{\begin{array}{cl}
    0, & \mbox{if} \,\, i = j \; \mbox{for each} \; i,j \in \an{L}, \; \mbox{for each} \; k \in \an{K}\\
    \frac{g^{k}_{ij}}{g^{k}_{ii}}, & \mbox{if} \,\, i \ne j  \; \mbox{for each} \; i,j \in \an{L}, \; \mbox{for each} \; k \in \an{K}
       \end{array}\right.
 \end{equation}
 and
 \begin{equation}\label{defbvdsl}
 \mathbf{v}= \displaystyle \left(\frac{n^{1}_1}{g^{1}_{11}}, \frac{n^{2}_1}{g^{2}_{11}},\dots, \frac{n^{K}_1}{g^{K}_{11}},
 \frac{n^{1}_2}{g^{1}_{22}}, \frac{n^{2}_2}{g^{2}_{22}},\dots,\frac{n^{K}_L}{g^{K}_{LL}} \right)\trans.
 \end{equation}

Note that $F$ is a block diagonal matrix under the assumption of synchronous transmission. More generally, asynchronous transmission can result in the $l$th user at tone $k$ having interference from the power allocation at neighboring tones of tone $k$. In this case, $F$ needs not be block diagonal.

As in the previous, the $l$th user is assigned a weight $w_l$ as an indicator of long-term priority. The problem of maximizing sum
 rate in a Gaussian interference-limited channel with $K$ frequency tones can be stated as the
 following maximum problem
\begin{equation}
\label{maxprobdsl0}
\max_{\sum_{k=1}^{K} p_l^{k} \le \bar{p}_l \; \forall \, l} \sum_{l=1}^L w_l \sum_{k=1}^K \log(1+\gamma_l^{k}(\p)),
\end{equation}
where $\gamma_l^{k}=g_{ll}^{k} p_l^{k}/(\sum_{j \ne l}g_{lj}^{k} p_j^{k}+n_l^{k})$. In order to be consistent with the previous results, where we assume $\mathbf{w}$ to be a probability vector, it is more convenient to consider the following equivalent problem:
\begin{equation}
\label{maxprobdsl}
\max_{\sum_{k=1}^{K} p_l^{k} \le \bar{p}_l \; \forall \, l} \sum_{l=1}^L \sum_{k=1}^K \tilde{w}_{lk} \log(1+\gamma_l^{k}(\p)).
\end{equation}
where $\tilde{w}_{lk}=(w_l/K)$ for all $l \in \an{L}, k \in \an{K}$, and $\mathbf{\tilde{w}} \in \R^{(K \times L)}$ is a probability vector.

As in Section \ref{sumratesection}, it is instrumental to consider an alternative formulation of
 (\ref{maxprobdsl}).
 \begin{theo}\label{altformpdsl}  The maximum problem
 (\ref{maxprobdsl}) is equivalent to the maximum problem.
 \begin{equation}
 \label{nonconvex1dsl}
 \begin{array}
 [c]{rl}
 \mbox{maximize} & \sum_l \sum_k \tilde{w}_{lk} \log(1+\gamma_l^{k})\\
 \mbox{subject to} &  \rho(\diag(\boldsymbol{\gamma})
 (F+(1/\bar{p}_l) \sum_k \mathbf{v}\mathbf{e}\trans_{(k \times l)})) \le 1 \,\;\; \forall \, l\in\an{L}, \\
 \mbox{variables:} & \gamma_l^{k}, \,\;\; \forall \, l, \; \forall \, k.
 \end{array}
 \end{equation}
 $\gam^{\star}$ is a maximal solution of the above problem if
 and only if $P(\gam^{\star})$ is a maximal solution
 $\p^{\star}$ of the problem (\ref{maxprobdsl}).  In particular,
 any maximal solution $\gam^{\star}$ satisfies the equality
 \begin{equation}
 \rho\left(\diag(\gam)\left(F+(1/\bar\p_l) \sum_k \v\e_{(k \times l)}\trans\right)\right) = 1
 \end{equation}
 for some integer $l\in [1,L]$.
 \end{theo}

Using Theorem \ref{altformpdsl}, it is straightforward to extend our three algorithms in Section \ref{algosection} for desigining optimal spectrum management schemes in DSL channels. We note however that our techniques are general enough to model asynchronous transmission, where the power allocations in neighboring tones interfere with one another. This results in intercarrier interference (ICI) and tone coupling in (\ref{ratedsl}); as mentioned earlier, ICI effects are captured by $F$ in (\ref{matrixFdsl}), where the optimal solution to (\ref{maxprobdsl}) is then characterized by the spectrum of different interference matrices.

In comparison to existing algorithms used for solving (\ref{maxprobdsl0}), we highlight briefly the advantages of our algorithms. Assuming synchronous transmission, the Optimal Spectrum Balancing (OSB) method is considered the state-of-the-art algorithm for computing the spectrum allocation \cite{Luo08}. The key idea of OSB is to leverage Lagrange dual decomposition to decouple the many uncoupled terms in (\ref{maxprobdsl0}) into $L$ lower dimensional subproblems (whose cost function on each tone is still nonconvex). However, this decoupling technique in the power domain cannot be done in the case of asynchronous transmission. On the other hand, in both synchronous and asynchronous transmission, Algorithm \ref{tayloralgo} and \ref{approxalgo} with the choice of supporting hyperplanes using the Friedland-Karlin inequalities permit decoupling of both users and frequency tones in the $\mathsf{SIR}$ domain.

\section{Conclusion}
We looked at optimizing power control to maximize the weighted sum rate in a Gaussian interference-limited channel that models multiuser communication in a CDMA wireless network or DSL cable binder. Using tools from irreducible nonnegative matrix theory, in particular the Perron-Frobenius Theorem and the Friedland-Karlin inequalities, we provided insights into the structural property of optimal power allocation strategies that maximize sum rates.
We showed that sum rate maximization can be restated as maximizing a convex function on a
 closed convex set. We proposed three algorithms to find the exact
 and approximate values of the optimal sum rates. In particular, our algorithms exploited the eigenspace of specially crafted nonnegative interference matrices, which, with the use of standard optimization algorithms, can provide useful upper bounds and feasible solutions to the nonconvex problem.

\appendix

 \section{Appendix: Friedland-Karlin results}

 In this section, we recall some results from
 \cite{Friedland75} and state the extensions of these
 results, and then illustrate their applications in this paper.  We first state the following
 extension of \cite[Theorem 3.1]{Friedland75}:
 \begin{theo}
 \label{boundonspectral}
 Let $A\in\R_+^{L\times L}$ be an irreducible matrix.
 Assume that $\x(A)=(x_1(A),\ldots,x_L(A))\trans,\y(A)=(y_1(A),\ldots,y_L(A))\trans>\0$
 are left and
 right Perron-Frobenius eigenvectors of $A$, normalized such that
 $\x(A)\circ\y(A)$ is a probability vector
 Suppose $\boldsymbol{\gamma}$ is a nonnegative vector. Then
 \begin{equation}
 \label{boundonspectraleqn}
 \rho(A)\prod_l \gamma_l^{(\mathbf{x}(A) \circ \mathbf{y}(A))_l}
 \le \rho(\diag(\gam)A).
 \end{equation}
 If $\gam$ is a positive vector then equality holds if and only if all
 $\gamma_l$ are equal.  Furthermore, for any positive vector
 $\z=(z_1,\ldots,z_L)\trans$, the following inequality holds:
 \begin{equation}\label{fkineq}
 \rho(A)\le \prod_{l=1}^L \left(\frac{(A\z)_l}{z_l}\right)^{(\mathbf{x}(A) \circ
 \mathbf{y}(A))_l}.
 \end{equation}
 If $A$ is an irreducible matrix with positive diagonal
 elements, then equality holds in (\ref{fkineq}) if and only if
 $\z=t\x(A)$ for some positive $t$.
 \end{theo}
 \proof
 Theorem 3.1 in \cite{Friedland75} makes the following
 assumptions.  First, in the inequality
 (\ref{boundonspectraleqn}), it assumes that $\gam>\0$.
 Second, in (\ref{fkineq}), it assumes that $\rho(A)=1$.
 Third, the equality case in (\ref{fkineq}) for $\z>\0$ is
 stated for a positive matrix $A$.
 We now show how to deduce the stronger version of Theorem 3.1
 claimed here.

 First, by using the continuity argument, we deduce
 the validity of (\ref{boundonspectraleqn}) for any $\gam\ge
 \0$.  Second,
 by replacing $A$ by $tA$, where $t>0$, we deduce that it
 is enough to show (\ref{fkineq}) in the case $\rho(A)=1$.

 Third, to deduce the equality case in (\ref{fkineq}) for
 $\z>\0$, we need to examine the proof of Lemma 3.2 in
 \cite{Friedland75}.  The proof of the Lemma 3.2 applies
 if the following condition holds.  For any sequence of probability vectors  $\z_i=(z_{1,i},
 \ldots,z_{L,i})\trans,
 i=1,\ldots$, which converges to a probability vector
 $\zet=(\zeta_1,\ldots,\zeta_L)\trans$, where $\zet$ has at least one zero
 coordinate, the function  $\prod_{l=1}^L \left(\frac{(A\z)_l}{z_l}\right)^{(\mathbf{x}(A) \circ
 \mathbf{y}(A))_l}$ tends to $\infty$ on the sequence $\z_i, i=1,\ldots$.  Assume that
 $\cA=\{l\in\an{L}, \zeta_l=0\}$.  Note that the complement of
 $\cA$ in $\an{L}$, denoted by $\cA^c$ is nonempty.

 Since $A=[a_{ij}]$ has positive diagonal entries, it follows
 that $\frac{(A\z)_l}{z_l}\ge a_{ll}>0$ for each $l\in \an{L}$.
 Since $A$ is irreducible, we deduce that there exist $l\in
 \cA$ and $m\in \cA^c$ such that $a_{lm}>0$.  Hence,
 $\lim_{i\to\infty}\frac{(A\z_i)_l}{z_{l,i}}=\infty$.  This shows
 that the unboundedness condition holds.

 \qed

 The following result gives an interpretation of the inequality
 (\ref{boundonspectraleqn}) in terms of the supporting
 hyperplane of the convex function $\log\rho(e^{\xit}B)$, where
 $B\in\R_+^{L\times L}$ is irreducible and $\xit\in\R^L$.
 \begin{theo}\label{suphyp} Let $B\in\R_+^{L\times L}$ be an irreducible matrix.
 Let $\et=(\eta_1,\ldots,\eta_L)\trans\in\R^L$
 satisfy the condition $\rho(e^{\et}B)=1$.  Denote $A=e^{\et}B$
 and assume that
 \par\noindent
 $\x(A)=(x_1(A),\ldots,x_L(A))\trans,\y(A)=(y_1(A),\ldots,y_L(A))\trans>\0$
 are left and
 right Perron-Frobenius eigenvectors of $A$, normalized such that
 $\x(A)\circ\y(A)$ is a probability vector.  Let
 \begin{equation}\label{defHxi}
 \rH(\xit)=\sum_{l=1}^L x_l(A)y_l(A)(\xi_l-\eta_l).
 \end{equation}
 Then $\rH(\xit)\le 0$ is the unique supporting hyperplane to the convex
 set $\log\rho(e^{\xit}B)\le 0$ at $\xit=\et$.

 \end{theo}
 \proof  Let $\xit\in\R^L$.  Then $e^{\xit}B=e^{\xit-\et}A$.
 Theorem \ref{boundonspectral} implies that $\rH(\xit)\le
 \log\rho(e^{\xit}B)$.  Thus, $\rH(\xit)\le 0$ if
 $\log\rho(e^{\xit}B)\le 0$.  Clearly, $\rH(\et)=0$.
 Hence, $\rH(\xit)\le 0$ is a supporting hyperplane of the
 convex set $\log\rho(e^{\xit}B)\le 0$.  Since the function
 $\log\rho(e^{\xit}B)$ is a smooth function of $\xit$, it follows that
 $\rH(\xit)\le 0$ is unique.  \qed

 We now give an application of (\ref{fkineq}) in Theorem \ref{boundonspectral}.
 \begin{theo}\label{apfkineq}  Let $B\in\R_+^{L\times L}$ be
 an irreducible matrix.  Let $\et=(\eta_1,\ldots,\eta_L)\trans\in\R^L$.
 Let $\w=\x(\diag(e^{\et}) B)\circ \y(\diag(e^{\et})
 B)=(w_1,\ldots,w_L)\trans$ be a probability vector.  Then for any positive vector
 $\z=(z_1,\ldots,z_L)\trans$
 \begin{equation}\label{apfkineq1}
 \sum_{l=1}^L w_l \log \frac{z_l}{(B\z)_l}\le
 -\log\rho(\diag(e^{\et})B) + \sum_{l=1}w_l\eta_l.
 \end{equation}
 If $B$ has a positive diagonal, then equality holds if and only
 if $\z=t\x(\diag(e^{\et})B)$ for some $t>0$.
 \end{theo}
 \proof Let $A=\diag(e^{\et})B$.  Then
 $$\sum_{l=1}^L w_l \log \frac{(B\z)_l}{z_l}=\sum_{l=1}^L w_l \log \frac{(A\z)_l}{z_l}
 - \sum_{l=1}w_l\eta_l.$$
 Use (\ref{fkineq}) to deduce (\ref{apfkineq1}).  The equality
 case follows from the equality case in (\ref{fkineq}).  \qed

 We now study the following inverse problem.
 \begin{prob}\label{invwprir}  Let $B\in\R_+^{L\times L}, \w\in\R_+^L$ be
 given irreducible matrix and positive probability vector, respectively.
 When does there exist $\et\in \R^L$ such
 that $\x(\diag(e^{\et})B)\circ \y(\diag(e^{\et})B)=\w$?
 If such $\et$ exists, when is it unique up to an addition
 $t\1$?
 \end{prob}

 To solve the inverse problem, we recall Theorem 3.2 in \cite{Friedland75}.
 \begin{theo}\label{thm3.2fk}  Let $A\in\R_+^{L\times L},
 \u,\v\in\R_+^{L}$ be given, where $A$ is irreducible with positive diagonal elements and
 $\u,\v$ are positive.  Then there exists
 $D_1,D_2\in\R_+^{L\times L}$ such that
 \begin{equation}\label{existD1D2}
 D_1AD_2\u=\u,\; \v\trans D_1AD_2=\v\trans,\;
 D_1=\diag(\f),\; D_2=\diag(\g) \textrm{ and } \f,\g
 >\0.
 \end{equation}
 The pair $(D_1,D_2)$ is unique to the change $(tD_1,
 t^{-1}D_2)$ for any $t>0$.  There exist $\et\in \R^L$ such
 that $\x(\diag(e^{\et}B))\circ \y(\diag(e^{\et}B))=\w$.
 Furthermore, $\et$ is unique up to an addition
 $t\1$.

 \end{theo}
 \begin{corol}\label{corthm3.2fk}  Let $B\in\R_+^{L\times L}, \w\in\R_+^L$ be
 given irreducible matrix with positive diagonal elements and positive probability vector,
 respectively.  Then there exists $\et\in \R^L$ such
 that $\x(\diag(e^{\et})B)\circ \y(\diag(e^{\et})B)=\w$.
 Furthermore, $\et$ is unique up to an addition of $t\1$.
 \end{corol}
 \proof  Let $\u=\1, \v=\w$.  Then there exists $D_1,D_2$ two
 diagonal matrices with positive diagonal entries such that
 $D_1BD_2 \1=\1, \w\trans D_1BD_2=\w\trans$.  Consider the
 matrix $D_2D_1 B=D_2(D_1 B D_2)D_2^{-1}$.  It is
 straightforward to see that $\x(D_2D_1 B)\circ
 \y(D_2D_1B)=\w$.  Hence, $\et$ is the unique solution of
 $\diag(e^{\et})=D_2D_1$.

 Assume that $\zet\in \R^L$ satisfies $\x(\diag(e^{\zet})B)\circ \y(\diag(e^{\zet})B)=\w$.
 By considering $\tilde\zet=\zet+t\1$, we may assume that
 $\rho(\diag(e^{\zet})B)=1$.
 Let $D_4=\diag(\x(\diag(e^{\zet})B)$.  Then
 $(D_4^{-1}\diag(e^{\zet})B D_4)\1=\1$.  Let
 $D_3=D_4^{-1}\diag(e^{\zet})$.  Hence, $\y(D_3BD_4)=\w$.  In
 view of Theorem \ref{thm3.2fk}, $\diag(e^{\zet})=D_4D_3=D_2D_1=\diag(e^{\et})$.  \qed

 Unfortunately, the matrix $F$ defined in (\ref{matrixF}) have
 zero diagonal entries and positive off-diagonal entries.
 For $L=2$, it is easy to show that
 $\x(F)\circ\y(F)=(\frac{1}{2},\frac{1}{2})\trans$.
 In particular, for $L=2$, Problem \ref{invwprir} is not
 solvable for $\w\ne (\frac{1}{2},\frac{1}{2})\trans$.
 Similarly, given positive $\u,\v\in \R^2$ such that
 $\u\circ\v\ne t(1,1)$ for any positive $t$, (\ref{existD1D2})
 does not hold for $A=F$.
 For $L\ge 3$, the situation is different.
 \begin{theo}\label{thm3.2ofdiagp}  Let $L\ge 3,A\in\R_+^{L\times L},
 \u=(u_1,\ldots,u_L)\trans,\v=(v_1,\ldots,v_L)\trans\in\R_+^{L}$ be given,
 where $A$ is a matrix with zero diagonal entries and positive
 off-diagonal elements, and
 $\u,\v$ are positive.  Assume that $\w=\u\circ\v$ is a
 probability vector satisfying the condition
 \begin{equation}\label{majcon}
 \sum_{\textrm{on all } j\ne l} w_j> w_l  \textrm{ for all } l\in\an{L}.
 \end{equation}
 Then there exists
 $D_1,D_2\in\R_+^{L\times L}$ such that (\ref{existD1D2})
 holds.
 \end{theo}
 \proof  Let $A_i=A+(1/i) I,i=1,\ldots$, where $I$ is
 the $L\times L$ identity matrix. Theorem \ref{thm3.2fk}
 implies
 \begin{eqnarray*}
 B_i=D_{1,i}A_iD_{2,i},\;
 D_{1,i}=\diag(\f_{i}),\;D_{2,i}=\diag(\g_{i}),\;
 B_i\u=\u,\;\v\trans B_i=\v\trans,\\
 \f_{i}=(f_{1,i},\ldots,f_{L,i})\trans,\;\w_{i}=(g_{1,i},\ldots,g_{L,i})\trans,\;
 s_i=\max_{j\in\an{L}} f_{j,i}=\max_{j\in\an{L}}
 g_{j,i},\;i=1,\ldots.
 \end{eqnarray*}
 Note that each entry of $B_i$ is bounded by $\frac{\max_j
 u_j}{\min_j u_j}$.
 By passing to the subsequence $B_{i_k},\f_{i_k},\g_{i_k}, 1\le
 i_1<i_2<\ldots$, we can assume that the first subsequence converges to $B$, and the
 last two subsequences converge in generalized sense:
 \begin{eqnarray*}
 \lim_{k\to\infty} B_{i_k}=B=[b_{jl}]\in \R_+^{L\times L},\;
 \lim_{k\to\infty} \f_{i_k}=\f=(f_1,\ldots,f_L)\trans,\;
 \lim_{k\to\infty} \g_{i_j}=\g=(g_1,\ldots,g_L)\trans, \\f_j,g_j\in
 [0,\infty],j=1,\ldots,L,
 \lim_{k\to\infty} s_{i_k}=s=\max_{j\in\an{L}} f_j=\max_{j\in\an{L}}
 g_j\in [0,\infty].
 \end{eqnarray*}
 Note that
 \begin{equation}\label{Buvid}
 B\u=\u,\quad \v\trans B=\v\trans.
 \end{equation}

 Assume first that $s<\infty$.  Then $B=\diag(\f)A\diag(\g)$.
 In view of (\ref{Buvid}), $\f \circ \g>\0$.  This proves the theorem in
 this case.

 Assume now that $s=\infty$.  Let
 \begin{eqnarray*}
 \cF_{\infty}=\{j\in\an{L},\;f_j=\infty\},\;
 \cF_{+}=\{j\in\an{L},\;f_j\in (0,\infty),\}\;
 \cF_{0}=\{j\in\an{L},\;f_j=0\},\\
 \cG_{\infty}=\{j\in\an{L},\;g_j=\infty\},\;
 \cG_{+}=\{j\in\an{L},\;g_j\in (0,\infty),\}\;
 \cG_{0}=\{j\in\an{L},\;g_j=0\}.
 \end{eqnarray*}
 Since off-diagonal entries of $A$ are positive, and $B\in\R_+^{L\times L}$
 it follows that $\cF_{\infty}=\cG_{\infty}=\{l\}$ for some
 $l\in\an{L}$.  Furthermore, $\cF_+=\cG_+=\emptyset$.  So
 $\cF_0=\cG_0=\an{L}\backslash\{l\}$.
 Assume first that $l=1$.  Then the principal submatrix
 $[b_{jl}]_{j=l=2}^L$ is zero.  (\ref{Buvid}) yields that
 $$b_{j1}=\frac{u_j}{u_1},\; b_{1j}=\frac{v_j}{v_1} \textrm{ for
 } j=2,\ldots,L,\; b_{11}u_1v_1+\sum_{j=2}^L u_jv_j =u_1v_1.
 $$
 Since $b_{11}\ge 0$, the above last equality contradicts the
 condition (\ref{majcon}) for $l=1$.  Similar argument implies the
 impossibility of  $\cF_{\infty}=\cG_{\infty}=\{l\}$ for any
 $l\ge 2$.  Hence, $s<\infty$ and we conclude the theorem.  \qed

 We do not know whether, under the conditions of Theorem \ref{thm3.2ofdiagp},
 the diagonal matrices $(D_1,D_2)$ are unique up to the
 transformation $(tD_1,t^{-1}D_2)$.
 We now generalize the above theorem.
 \begin{theo}\label{thm3.2ofdiagp1}  Let
 $$L\ge 2,\; A=[a_{jl}]_{j=l=1}^L\in\R_+^{L\times L},\;
 \0<\u=(u_1,\ldots,u_L)\trans,\;\v=(v_1,\ldots,v_L)\trans\in\R_+^{L}$$ be
 given.  Assume that $A$ has positive off-diagonal elements, and $\w=\u\circ\v$ is a
 probability vector satisfying the condition
 \begin{equation}\label{majcon1}
 \sum_{\textrm{on all } j\ne l} w_j> w_l
 \end{equation}
 for each $l$  such that $a_{ll}=0$.
 Then there exists
 $D_1,D_2\in\R_+^{L\times L}$ such that (\ref{existD1D2})
 holds.
 \end{theo}
 \proof
 Assume first that $L\ge 3$.
 In view of Theorems \ref{thm3.2fk} and
 \ref{thm3.2ofdiagp},
 it is enough to assume that $A$ has positive and zero diagonal
 entries.  Apply the proof of Theorem \ref{thm3.2ofdiagp} and
 the following observation.
 If  $\cF_{\infty}=\cG_{\infty}=\{l\}$ then
 $a_{ll}=0$.

 Assume now that $L=2$.  Note that if $A$ has a zero diagonal
 then the condition (\ref{majcon}) can not hold.  Assume now
 that $A$ has at least one positive diagonal element.  Then the
 above arguments for $L\ge 3$ apply.
 \qed
\begin{corol}\label{corthm3.2ofdiag}  Let $B=[b_{jl}]_{j=l=1}^L\in\R_+^{L\times L}, \w\in\R_+^L$ be
 given matrix with positive off-diagonal elements and a positive probability vector,
 respectively.  Assume that $L\ge 2$ and $\w$ satisfies the conditions
 (\ref{majcon1}) for each $l$ such that $b_{ll}=0$.
 Then there exists $\et\in \R^L$ such
 that $\x(\diag(e^{\et})B)\circ \y(\diag(e^{\et})B)=\w$.
 \end{corol}

 \end{document}